# Kenmotsu Manifolds Admitting a Non-Symmetric Non-Metric Connection

S. K. Chaubey[1], S. K. Yadav[2] and *Mahesh Garvandha*[3]

[1, 3] Mathematics Section, Department of Information Technology, Shinas College of Technology, PO Box 77, Postal Code 324 Al Aqur, Shinas, Sultanate of Oman.

[2] Department of Mathematics, Poornima college of Engineering, Rajasthan, India.

e-mail[1] : sudhakar.chaubey@shct.edu.om;   e-mail[2] : prof_sky16@yahoo.com;

e-mail[3] : mahesh.garvandha@shct.edu.om

**Abstract**

The aim of the present paper is to study the properties of Kenmotsu manifolds equipped with a non-symmetric non-metric connection. We also establish some curvature properties of Kenmotsu manifolds. It is proved that a Kenmotsu manifold endowed with a non-symmetric non-metric is irregular.

1. **Literature Review**

A contact metric manifold is capable to resolve many issues of sciences, engineering and medical sciences, and hence it is attracting the researchers to work in this area. Boothby & Wang (1958) started the study of a differentiable manifold with contact and almost contact metric structures. Kenmotsu (1971) introduced a class of almost contact metric manifold and named as Kenmotsu manifold. Since then, the properties of Kenmotsu manifolds have studied by several authors such as De & Pathak (2004), Sinha & Srivastava (1991), Jun, De & Pathak (2005), De, Yildiz and Yaliniz (2008), Chaubey et al. (2010, 2012, 2015, 2018), De (2008), Cihan (2006) and many others.

Let $M$ be a Riemannian manifold associated with the Riemannian metric $g$. A linear connection $\tilde{\nabla}$ on $M$ is said to be symmetric if the torsion tensor $\tilde{T}$ of $\tilde{\nabla}$ vanishes, otherwise it is non-symmetric.

If the torsion tensor $\tilde{T}$ assumes the form $\tilde{T}(X,Y) = \pi(Y)X - \pi(X)Y$ for all vector fields $X$ and $Y$ on $M$, then then linear connection $\tilde{\nabla}$ is called semi-symmetric connection (Friedmann, & J. A. Schouten, 1924). Moreover, if $\tilde{\nabla} g = 0$ on $M$ then $\tilde{\nabla}$ is said to be metric, otherwise non-metric (Agashe & Chafle, 1992). Chaubey & Ojha (2008) defined and studied the properties of non-symmetric non-metric connection on almost contact metric manifolds. The geometrical properties of the same connection have been studied by several authors. We refer (Chaubey & Kumar (2010); Chaubey & De (2019); Chaubey, & Yildiz (2019); Chaubey et al. (2019)) and their references.

Above studies motivate us to study the properties of Kenmotsu manifold equipped with a non-symmetric non-metric connection.

## 2. Preliminaries:

A smooth manifold $M$ of dimension (2n+1) is said to be an almost contact metric manifold if it admits a (1, 1) tensor field $\phi$, (1, 0) type vector field $\xi$, (0, 1) type vector field $\eta$ and a compatible metric $g$ of type (0, 2) satisfies

$$\phi^2 = -I + \eta \otimes \xi, \ \eta(\xi) = 1 \ \text{and} \ g(\phi X, \phi Y) = g(X,Y) - \eta(X)\eta(Y) \qquad (2.1)$$

for all $X$ and $Y$ on $M$ (Blair ,1976). Additionally, if $M$ satisfies

$$\nabla_X \xi = X - \eta(X)\xi \Leftrightarrow (\nabla_X \eta)(Y) = g(X,Y) - \eta(X)\eta(Y) \qquad (2.2)$$

for all $X$ on $M$, then $M$ is said to be a Kenmotsu manifold (1971). Here $\nabla$ denotes the Levi-Civita connection of $g$. It is observed that the manifold holds the following relations,

$$\eta(R(X,Y)Z) = \eta(Y)g(X,Z) - \eta(X)g(Y,Z), \qquad (2.3)$$

$$R(X,Y)\xi = \eta(X)Y - \eta(Y)X, \qquad (2.4)$$

$$R(\xi,X)Y = \eta(Y)X - g(X,Y)\xi, \qquad (2.5)$$

$$S(X,\xi) = -2n \ \eta(X) \qquad (2.6)$$

for all $X$, $Y$ and $Z$ on $M$ (Kenmotsu (1971), Chaubey & Ojha (2010), Chaubey & R. H. Ojha (2012)). Here $R$ and $S$ denote the Riemannian curvature and Ricci tensors of $g$, respectively.

A Kenmotsu manifold $M$ is said to be $\eta$-Einstein if the non-vanishing Ricci tensor $S$ satisfies the relation $S(X,Y) = ag(X,Y) + b\eta(X)\eta(Y)$ for all $X$ and $Y$ on $M$, where $a$ and $b$ are smooth functions on $M$ (Kenmotsu (1971)).

### 3. Non-symmetric non-metric connection

Let $M$ be a Kenmotsu manifold of dimension $(2n+1)$. A linear connection $\tilde{\nabla}$ on $M$, defined by

$$\tilde{\nabla}_X Y = \nabla_X Y - \eta(Y)X - g(X,Y)\xi \tag{3.1}$$

for all vector fields $X$ and $Y$ on $M$, is known as a non-symmetric non-metric connection (De & Pathak (2004)), if the torsion tensor $\tilde{T}$ of $\tilde{\nabla}$ takes the form $\tilde{T}(X,Y) = \pi(X)Y - \pi(Y)X$ and

$$(\tilde{\nabla}_X g)(Y,Z) = 2\eta(Y)g(X,Z) + 2\eta(Z)g(X,Y). \tag{3.2}$$

In consequence of the equations (2.1), (2.3) and (3.1), we get

$$\tilde{\nabla}_X \xi = \nabla_X \xi - X - \eta(X)\xi = -2\eta(X)\xi. \tag{3.3}$$

If $\tilde{R}$ denotes the Riemannian curvature tensor with respect to $\tilde{\nabla}$, then it relates to $R$ by the relation

$$\tilde{R}(X,Y)Z = R(X,Y)Z - \beta(X,Z)Y + \beta(Y,Z)X - g(Y,Z)(\nabla_X \xi - \eta(X)\xi)$$

$$+ g(X,Z)(\nabla_Y \xi - \eta(Y)\xi), \tag{3.4}$$

where $\beta$ is tensor field of type $(0,2)$ and defined as

$$\beta(X,Y) = (\nabla_X \eta)(Y) + \eta(X)\eta(Y) + g(X,Y) = 2g(X,Y). \tag{3.5}$$

In view of (3.5), equation (3.4) takes the form

$$\tilde{R}(X,Y)Z = R(X,Y)Z + g(Y,Z)X - g(X,Z)Y + 2[g(Y,Z)\eta(X) - g(X,Z)\eta(Y)]\xi. \tag{3.6}$$

Contracting equation (3.6) along the vector field $X$, we lead

$$\tilde{S}(Y,Z) = S(Y,Z) + 2(n+1)g(Y,Z) - 2\eta(Y)\eta(Z), \quad (3.7)$$

which gives

$$\tilde{Q}Y = QY + 2(n+1)Y - 2\eta(Y)\xi. \quad (3.8)$$

The contraction of (3.8) gives

$$\tilde{r} = r + 2n(2n+3). \quad (3.9)$$

Here $\tilde{r}$ and $r$ are the scalar curvatures with respect to $\tilde{\nabla}$ and $\nabla$, respectively. $\tilde{Q}$ and $Q$ denote the Ricci operators corresponding to the Ricci tensors $\tilde{S}$ and $S$ with respect to $\tilde{\nabla}$ and $\nabla$, respectively.

Setting $Z = \xi$ in (3.7) and then using the equations (2.1) and (2.4), we obtain

$$\tilde{R}(X,Y)\xi = R(X,Y)\xi + \eta(Y)X - \eta(X)Y + 2[\eta(Y)\eta(X) - \eta(X)\eta(Y)]\xi = 0.$$

This shows that the Kenmotsu manifold M is irregular ($\tilde{R}(X,Y)\xi = 0$) with respect to $\tilde{\nabla}$. Thus we state:

**Theorem 3.1.** *Every $(2n+1)$-dimensional Kenmotsu manifold equipped with $\tilde{\nabla}$ is irregular with respect to $\tilde{\nabla}$.*

4. **Ricci semi-symmetric Kenmotsu manifold with a non-symmetric non-metric connection**

This section deals with the study of Ricci semi-symmetric Kenmotsu manifold equipped with a non-symmetric non-metric connection $\tilde{\nabla}$. It is well known that

$$(\tilde{R}(X,Y) \cdot \tilde{S})(Z,U) = -\tilde{S}(\tilde{R}(X,Y)Z,U) - \tilde{S}(Z,\tilde{R}(X,Y)U).$$

In view of (2.1) and (3.7), we obtain

$$(\tilde{R}(X,Y) \cdot \tilde{S})(Z,U) = (R(X,Y) \cdot S)(Z,U) - g(Y,Z)S(X,U) + g(X,Z)S(Y,U) -$$

$$g(Y,U)S(Z,X) + g(X,U)S(Z,Y), \tag{4.1}$$

where $(R(X,Y) \cdot S)(Z,U) = -S(R(X,Y)Z,U) - S(Z, R(X,Y)U)$ for all vector fields X, Y, Z and U on M. If possible, we suppose that $\tilde{R} \cdot \tilde{S} = R \cdot S$, then (4.1) gives

$$-g(Y,Z)S(X,U) + g(X,Z)S(Y,U) - g(Y,U)S(Z,X) + g(X,U)S(Z,Y) = 0 \tag{4.2}$$

Changing Y and U by $\xi$ in (4.2) and then using the equations (2.1), (2.2), and (2.6), we obtain

$$S(Z,X) = -2ng(Z,X), \quad r = -2n(2n+1), \tag{4.3}$$

which shows that the Kenmotsu manifold endowed with a non-symmetric non-metric connection $\tilde{\nabla}$, under assumption is an Einstein manifold. Also, from equations (3.7) and (4.3), we find

$$\tilde{S}(Z,X) = 2g(Z,X) - 2\eta(Z)\eta(X). \tag{4.4}$$

This reflects that the Kenmotsu manifold M w.r.t. $\tilde{\nabla}$ under consideration is an $\eta$-Einstein manifold. Thus, we can state:

**Theorem 4.1.** *Let $(M,g)$ be a $(2n+1)$-dimensional Kenmotsu manifold equipped with a non-symmetric non-metric connection $\tilde{\nabla}$. If $\tilde{R} \cdot \tilde{S} = R \cdot S$ on M, then the manifold to be an Einstein manifold also, M is an $\eta$-Einstein with respect to $\tilde{\nabla}$.*

With the help of (4.3), equation (3.9) takes the form $\tilde{r} = 4n$, which shows that if a Kenmotsu manifold equipped with $\tilde{\nabla}$ satisfies $\tilde{R} \cdot \tilde{S} = R \cdot S$, then the scalar curvature with respect to $\tilde{\nabla}$ is constant. Thus, we state:

**Corollary.** *If a $(2n+1)$-dimensional Kenmotsu manifold M endowed with $\tilde{\nabla}$ satisfies $\tilde{R} \cdot \tilde{S} = R \cdot S$, then the scalar curvature of M with respect to $\tilde{\nabla}$ to be constant.*

Arslan et al. (2014) proved that

"A semi-Riemannian Einstein manifold M of dimension $n, \geq 4$, satisfies

$$R \cdot C - C \cdot R = \frac{r}{n(n-1)} Q(g, R) = \frac{r}{n(n-1)} Q(g, C), \qquad (4.5)$$

where $Q(g, R)$ denotes the *Tachibana tensor* and $C$ is the *conformal curvature tensor*."

From (4.3) and (4.5), we have

$$C \cdot R - R \cdot C = Q(g, R) = Q(g, C). \qquad (4.6)$$

Thus, we state:

**Corollary.** *If a* $(2n + 1)$-*dimensional Kenmotsu manifold M equipped with* $\tilde{\nabla}$ *satisfies*

$\tilde{R} \cdot \tilde{S} = R \cdot S$ , *then M satisfies the equation (4.6)*.